\documentclass[twocolumn,10pt]{asme2e}
\usepackage{amsmath,amssymb,amstext}
\usepackage{graphicx,color}
\usepackage{graphics}
\confshortname{GT2013} \conffullname{ASME Turbo Expo 2013}

\confdate{June 3-7} \confyear{2013} \confcity{San Antonio, Texas}
\confcountry{USA}

\title{RISK ESTIMATION FOR LCF CRACK INITIATION}
\author{Sebastian Schmitz\thanks{Address all correspondence to this author. Second affiliation: Institute of
Computational Science, Universita della Svizerra Italiana, Lugano,
Ticino, 6900, Switzerland, Email: sebastian.schmitz@usi.ch}
    \affiliation{
    Gas Turbine Department of Materials and Technology\\
    Siemens AG Energy\\
    M\"ulheim an der Ruhr, Nordrhein-Westfalen, 45473\\
    Germany\\
    Email: schmitz.sebastian@siemens.com\vspace{5mm}
    }\\
    {\tensfb Georg Rollmann}
    \affiliation{Gas Turbine Department of Materials and Technology\\
    Siemens AG Energy\\
    M\"ulheim an der Ruhr, Nordrhein-Westfalen, 45473\\
    Germany\\
    Email: georg.rollmann@siemens.com
    }
}

\author{Hanno Gottschalk%\thanks{Address all correspondence to this author.} \\
\affiliation{Faculty of Mathematics and Natural Science\\
    Bergische Universit\"at Wuppertal\\
    Wuppertal, Nordrhein-Westfalen, 42097\\
    Germany\\
    Email: hanno.gottschalk@uni-wuppertal.de\vspace{5mm}
    }\\
       {\tensfb Rolf Krause}
    \affiliation{Institute of Computational Science\\
    Universita della Svizzerra Italiana\\
    Lugano, Ticino, 6900\\
    Switzerland\\
    Email: rolf.krause@usi.ch
    }
}

\begin{document}

\maketitle

%%%%%%%%%%%%%%%%%%%%%%%%%%%%%%%%%%%%%%%%%%%%%%%%%%%%%%%%%%%%%%%%%%%%%%
\begin{abstract}
{\it An accurate risk assessment for fatigue damage is of vital
importance for the design and service of today's turbomachinery
components. We present an approach for quantifying the probability
of crack initiation due to surface driven low-cycle fatigue (LCF).
This approach is based on the theory of failure-time processes and
takes inhomogeneous stress fields and size effects into account.
The method has been implemented as a finite-element postprocessor
which uses quadrature formulae of higher order. Results of
applying this new approach to an example case of a gas-turbine
compressor disk are discussed.}
\end{abstract}

%%%%%%%%%%%%%%%%%%%%%%%%%%%%%%%%%%%%%%%%%%%%%%%%%%%%%%%%%%%%%%%%%%%%%%
%\begin{nomenclature}
%\entry{A}{You may include nomenclature here.}
%\entry{$\alpha$}{There are two arguments for each entry of the
%nomemclature environment, the symbol and the definition.}
%\end{nomenclature}

%The spacing between abstract and the text heading is two line
%spaces. The primary text heading is  boldface in all capitals,
%flushed left with the left margin.  The spacing between the  text
%and the heading is also two line spaces.

%%%%%%%%%%%%%%%%%%%%%%%%%%%%%%%%%%%%%%%%%%%%%%%%%%%%%%%%%%%%%%%%%%%%%%
\section*{INTRODUCTION}

Due to the necessity for a more flexible service of gas turbines,
low-cycle fatigue (LCF) design has become of essential importance
in today's gas turbine engineering. The end of LCF-life of an
engineering component is often defined by the initiation of a
crack of a certain size. LCF is mostly surface driven so that an
LCF crack initiates on the component's surface. In the case of
polycrystalline metal the grain structure has a great influence on
the LCF failure mechanism which is of stochastic nature. This can
result in a statistical scatter of a factor of 10 between the
highest and lowest load cycles to crack initiation, even under lab
conditions.

%given a temperature and strain level.
%Stochastic attributes prevail on the macroscopic scale, where a
%statistical scatter of a factor 10 between the highest and lowest
%load cycles to crack initiation is a common phenomenon, even under
%lab conditions.

Standard design approaches often derive a predicted component life
with respect to LCF from the average life times of the most loaded
points on the component plus safety factors which account for the
scatter band, size effects and uncertainties in the stress and
temperature fields. However, this method can have a lack of
sufficient accuracy which can result in designs that are too
conservative or too optimistic, for example. In particular, the
Coffin-Manson-Basquin equation and W\"ohler curves play an
essential role in reliability estimations regarding fatigue.
%confer \cite{Harders_Roesler}.
Moreover, size effects have a
significant influence on fatigue life and are mostly considered by
safety factors.
%For example, different geometries of the same
%material and under the same temperature and stress fields can lead
%to different W\"ohler curves, confer \cite{Vormwald}. %Section 3.5
For more detailed discussions of fatigue
%in general, applications
%in engineering and on failure mechanism with respect to different
%material classes
confer \cite{Harders_Roesler}, \cite{Vormwald},
\cite{Fedelich} and \cite{Sornette}.

In this work we focus on LCF in conjunction with polycrystalline
metals and in particular on the number of load cycles until crack
initiation. We present a local and probabilistic model for LCF
according to \cite{Gottschalk_Schmitz} and
\cite{Siemens_Juelich_Wuppertal_Lugano} which we use to estimate
the risk for LCF crack initiation on a compressor disk.
%which we have implemented as a postprocessor for finite element
%analysis (FEA).
Here, the reaction of a component to cyclic loads is taken into
account via a linear elastic finite element analysis (FEA) and via
Neuber shakedown, see \cite{Elasticity_Ciarlet_1},
\cite{Finite_Elements_Ern} and \cite{Knop_Jones_Molent_Wang},
respectively.

%Metallic material usually does not consist of one single crystal
%but of different crystalline regions which are called grains,
%confer \cite{Harders_Roesler}. After the melt cooling and
%solidification of metal has resulted in recrystallization, small
%solidified regions with crystalline structure form.  Since the
%crystal nuclei are independent they have no long-range order in
%common. This is an important characteristic of polycrystalline
%metal which consists of grains with an order of magnitude
%typically in the micrometer and millimeter scale.
%%The crystal boundaries between the grains are regions where
%%different crystalline orientations abut.

Because LCF cracks are small in the initiation phase they will
only influence stress fields on micro- and mesoscales and so it is
reasonable to assume that crack formation in one region of the
component's surface is not influenced by the crack forming process
on another part. Thus crack formation can be considered as a
problem of spatial statistics \cite{Sherman}. Then,
\cite{Gottschalk_Schmitz} and
\cite{Siemens_Juelich_Wuppertal_Lugano} infer that hazard rates
for crack initiation have to be integrals over some local function
depending on the local stress or strain fields. Therefore, hazard
rates for the component can be expressed by a surface integral
over some crack formation intensity function depending on local
fields. The latter can also be regarded as the density for the
intensity measure of a Poisson point process (PPP) in time and
space, confer \cite{Baddeley} and \cite{Klenke}. In
\cite{Fedelich} the role of the PPP is also emphasized in that
context.

Following \cite{Gottschalk_Schmitz} and
\cite{Siemens_Juelich_Wuppertal_Lugano} we model the hazard rates
with a rather conservative approach by means of intensity measures
of
Weibull type. Note that %location-scale-based
scale-shape distributions are very common in reliability
statistics, confer \cite{Escobar_Meeker}. Furthermore, we assume
that the scale variable $N_{\textrm{det}}$ is of the same
functional form as the usual Coffin-Manson-Basquin equation,
confer \cite{Gottschalk_Schmitz},\cite{Harders_Roesler} and
(\ref{materials_1.1}) below. This results in a Weibull
distribution for the number $N$ of cycles of first crack
initiation on the component as well.

In contrast to \cite{Fedelich} this purely phenomenological
approach avoids detailed modeling at the meso scale which
facilitates calibration with experiments. Both approaches have the
use of Weibull distributions in common. From a materials
engineering point of view the model of \cite{Gottschalk_Schmitz}
and \cite{Siemens_Juelich_Wuppertal_Lugano} has the significant
advantage -- compared to standard methods in fatigue -- of
bypassing the standard specimen approach and considering size
effects. Besides LCF-test results with standardized specimens,
different strain-controlled LCF results can be used together for
the calibration of the model such as results from specimens with
different geometries or under different inhomogeneous strain
fields.
%Standard specimens of course can still be used to calibrate the
%model. But other strain-controlled LCF results could be used as
%well for that purpose such as results from specimens under
%inhomogeneous stress fields.

In order to numerically compute the Weibull distribution,
quadrature formulae are employed for the corresponding
integration. As locations of stress concentrations result in
higher nonlinearities in the integrand we use quadrature formulae
of higher orders.
%Moreover, the method of reduced integration is
%often used in FEA packages where certain integration points are
%neglected to reduce needed computing performance. Thus, we will
%not apply integration points of FEA packages.
For this purpose we interpolate the field values according to
principles of mesh and finite-element generation, confer
\cite{Finite_Elements_Ern}, \cite{Finite_Elements_Ciarlet_1} and
\cite{Braess}.
%Here, geometric and functional transformations from the reference
%finite element to the finite elements of the component play an
%essential role.

Having obtained the Weibull shape and scale parameter from
calibration and from numerical integration, respectively, the
corresponding distribution function yields the probability for LCF
crack initiation with respect to the number $N$ of load cycles. In
the design process one can decide which number $N$ is acceptable
corresponding to the assigned risk. Another important quantity of
the LCF model is the local crack initiation density. This field
shows how much each region of the surface contributes to the
overall expected number of crack initiations and leads to critical
as well as to possibly overengineered parts of the component.

This paper combines some aspects of materials engineering,
reliability statistics and FEA. In the first section we consider
linear elasticity, Neuber shakedown, fatigue analysis and we
present the local and probabilistic model for LCF. In Section
\ref{section_Postprocessor} we first focus on FEA. Then, we
discuss numerical integration and important functions of our
presented approach. The last section shows results of applying our
method to an example case of a gas-turbine compressor disk.

%%%%%%%%%%%%%%%%%%%%%%%%%%%%%%%%%%%%%%%%%%%%%%%%%%%%%%%%%%%%%%%%%%%%%%
\section{A LOCAL AND PROBABILISTIC MODEL FOR LCF}

In this section we discuss results of linear isotropic elasticity,
Neuber shakedown and fatigue analysis. The last subsection
presents the local and probabilistic model for LCF which is
introduced and motivated in \cite{Gottschalk_Schmitz} and
\cite{Siemens_Juelich_Wuppertal_Lugano}.

%%%%%%%%%%%%%%%%%%%%%%%%%%%%%%%%%%%%%%%%%%%%%%%%%%%%%%%%%%%%%%%%%%%%%%
\subsection{Linear Isotropic Elasticity and
Neuber Shakedown}\label{section_thermoelasticity_plasticity}

In this work we consider designs made from single-phased
polycrystalline metal. Usually metallic material does not consist
of one single crystal but of different crystalline regions which
are called grains. These grains have an order of magnitude
typically in the micrometer and millimeter scale. The continuum
mechanical approach assumes that the considered sizing scale is
large compared to the inter-atom distances. Thus the material is
considered to be smeared and all quantities are continuous. In
general, single crystals have anisotropic properties. In
polycrystals the orientation of the grains is randomly
distributed. If the grains are sufficiently small compared to the
component's size, the anisotropic effects of the grains average
out and an approximately isotropic material can be assumed.

In the following, we employ the continuum mechanical approach and
assume isotropic material behavior at scales significantly larger
than the grain size and assume sufficiently small deformations.
Thus, linear isotropic elasticity can be applied to describe the
behavior of components from single-phased polycrystalline metal
under external loading and plasticity can be considered by Neuber
shakedown. In this section we present theoretical backgrounds and
closely follow Section 2 of \cite{Gottschalk_Schmitz} which is
based on \cite{Hetnarski} and \cite{Harders_Roesler}. Note that
our example case of a compressor disk is subject to a homogeneous
temperature field so that we do not consider thermoelasticity.

Let $\Omega$ be a domain which represents the component shape
filled with a deformable medium such as polycrystalline metal
which is initially at equilibrium. Moreover, let $\nu$ be the
normal on the surface $\partial\Omega$ of $\Omega$, let
$\textbf{f}$ be an external load and let $\textbf{u}$ be the
three-dimensional displacement field in $\Omega$. Finally, let
$\partial\Omega_D,\partial\Omega_N$ be a partition of the boundary
where $\partial\Omega_D$ is clamped and on $\partial\Omega_N$ a
normal load $\textbf{g}$ is imposed. Then, according to
\cite{Finite_Elements_Ern} the mixed boundary value problem (BVP)
of linear isotropic elasticity is described by:
\begin{equation}\label{Thermoelasticity.1.0}
\nabla\cdot\sigma^e(\textbf{u})+\textbf{f}=0\quad\textrm{in
}\Omega
\end{equation}
with
$\sigma^e(\textbf{u})=\lambda(\nabla\cdot\textbf{u})I+\mu(\nabla
\textbf{u}+\nabla \textbf{u}^T)$ and with boundary conditions
$\textbf{u}=0$ on $\partial\Omega_D$ and
$\sigma^e(\textbf{u})\cdot \nu=\textbf{g}$ on $\partial\Omega_N$.
Here, $\lambda$ and $\mu$ are the Lame coefficients. The
linearized strain rate tensor $\varepsilon^e(\textbf{u})$ is
defined as $\varepsilon^e(\textbf{u})=\frac{1}{2}(\nabla
\textbf{u}+\nabla \textbf{u}^T)$, i.e.
$\varepsilon_{ij}^e=\frac{1}{2}\left(\frac{\partial
\textbf{u}_i}{\partial x_j}+\frac{\partial \textbf{u}_j}{\partial
x_i}\right)$ for $i,j=1,2,3$. Numerical solutions of the BVP can
be computed by an FEA, confer \cite{Hetnarski},
\cite{Finite_Elements_Ern} and Section \ref{section_Postprocessor}
below.

%Moreover, confer Section \ref{section_Postprocessor} where we
%address important concepts of FEA and describe our postprocessor
%namenlos which employs results from an FEA model of a component.

%In the following we phenomenologically discuss time-independent
%plasticity.
The knowledge of the threshold between elastic and plastic
deformations is very important as plastic deformations can allude
to an imminent residual fracture. According to
\cite{Harders_Roesler} this threshold is often described by
so-called yield criteria.
%For multi-axile stresses scalar comparison stresses $\sigma_V$ are
%introduced which depends on the stress tensor. If the comparison
%stress reaches a critical value $\sigma_{\textrm{crit}}$ the
%material begins to yield. A yield criteria can have the form
%$y(\sigma_{ij})=0$, where
%$y(\sigma_{ij})=\sigma_V(\sigma_{ij})-\sigma_{\textrm{crit}}$.
In this work we use the von Mises yield criterion which is given
by
\begin{equation}\label{Thermoelasticity.1.1}
\sqrt{\frac{1}{6}\left[(\sigma_{1}-\sigma_{2})^2+(\sigma_{1}-\sigma_{3})^2+(\sigma_{2}-\sigma_{3})^2\right]}=k_F,
\end{equation}
where $k_F$ is the critical value of the criterion and
$\sigma_1,\sigma_2,\sigma_3$ are the principal stresses.
%where $\sigma:\sigma=\sigma^{ij}\sigma_{ij}$ with Einstein's convention of summation.
%given in the space of principal stresses.
Here, the left-hand side is proportional to the elastic strain
energy of distortion. If the criterion is applied to uniaxial
tensile tests the relationship $k_F=R_p/\sqrt{3}$ is obtained,
where $R_p$ is the critical value of the only nonzero principle
stress in a uniaxial tensile test. The von Mises stress is defined
as
\begin{equation}\label{Thermoelasticity.7.1}
%\begin{split}
\sigma_v=\sqrt{\frac{1}{2}\left[(\sigma_{1}-\sigma_{2})^2+(\sigma_{1}-\sigma_{3})^2+(\sigma_{2}-\sigma_{3})^2\right]}.
%=\sqrt{\frac{3}{2}\sum_{i,j=1}^3\sigma'_{ij}\sigma'_{ij}}
%\end{split}
\end{equation}
Thus, the previous criterion can be written as $\sigma_v=R_p$ and
be used to predict yielding of metal under any loading condition
from results of uniaxial tensile tests.

In the following we introduce the Ramberg-Osgood equation, confer
\cite{Ramberg}. It can be used to locally derive strain levels
from scalar comparison stresses such as the von Mises stress.
These strain levels determine strain-controlled fatigue life. The
Ramberg-Osgood equation establishes stress-strain curves of metals
near their yield points. It is very accurate in the case of smooth
elastic-plastic transitions which can be observed for metals that
harden with plastic deformations, for example. If $K$ denotes the
cyclic strain hardening coefficient and $n$ the cyclic strain
hardening exponent the Ramberg-Osgood equation is given by
\begin{equation}\label{Thermoelasticity.8.1}
\varepsilon_v=\frac{\sigma_v}{E}+\left(\frac{\sigma_v}{K}\right)^{1/n}
\end{equation}
with Young' modulus $E=\frac{\mu(3\lambda+2\mu)}{\lambda+\mu}$.
The equation defines the comparison strain $\varepsilon_v$, where
we also write $\varepsilon_v=RO(\sigma_v)$.

Finally, we present the method of Neuber shakedown, confer
\cite{Neuber} and \cite{Harders_Roesler}. Let $\sigma^e_v$ denote
the von Mises stress which is obtained only from linear elastic
computations and let $\sigma_v$ be the von Mises stress which also
considers plasticity. Here, $\sigma_v$ is called elastic-plastic
von Mises stress. If linear elasticity leads to stress values
greater than the material yield strength, Neuber shakedown will
estimate corresponding elastic-plastic stress values. An
energy-conservation ansatz is the foundation of Neuber's approach.
It results in a relationship between the elastic von Mises
stress\footnote{Confer \cite{Hoffmann} for details on Neuber
shakedown in conjunction with equivalent stresses. As an
alternative to this method one could also use Glinka's method,
confer e.g. \cite{Knop_Jones_Molent_Wang}.} $\sigma^e_v$ and the
elastic-plastic von Mises stress $\sigma_v$:
\begin{equation}\label{Thermoelasticity.9.1}
\frac{(K_t\,\sigma^e_v)^2}{E}=\sigma_v\,\varepsilon_v=\frac{\sigma_v^2}{E}+\sigma_v\left(\frac{\sigma_v}{K}\right)^{1/n}.
\end{equation}
Here, the Ramberg-Osgood approach is also used. $K_t$ is the notch
factor, which is set to one if $\sigma^e_v$ is obtained from the
BVP (\ref{Thermoelasticity.1.0}) where notches are incorporated in
the boundary definition. Given the elastic comparison stress
$\sigma_v^e$, we can thus calculate the elastic-plastic von Mises
stress by solving (\ref{Thermoelasticity.9.1}). Thereby, we are
able to obtain $\varepsilon_v$ from (\ref{Thermoelasticity.8.1}).
Note that we also write $\sigma_v=SD^{-1}(\sigma_v^e)$.

%For the case of a very fine structure of grains in polycrystalline
%metal the property is supposed to be satisfied well

%\begin{figure}
%\begin{minipage}[hbt]{8cm}
%    \centering
%\scalebox{0.3}{\input{Doktorarbeit_hysteresis_2.pstex_t}}
%%\vspace{0.25cm}
%\caption{Hysteresis and cyclic loading.}
%    \label{Figure_hysteresis}
%\end{minipage}
%\end{figure}

\subsection{Fatigue and the Coffin-Manson-Basquin Equation}\label{section_fatigue}

Fatigue describes the damage or failure of material under cyclic
loading, confer \cite{Harders_Roesler},\cite{Schott} and
\cite{Vormwald}.
%One reason for the importance of fatigue analysis
%is that material can be damaged by much lower load amplitudes of
%cyclic loading compared to the static case. Another reason are
%less significant plastic deformations in a ductile material under
%cyclic loading before failure.
%Thus, an imminent damage is more difficult to detect in that case.
Major examples of cases where fatigue occurs are activation and
deactivation operations, e.g. of motor vehicles and of gas
turbines, and oscillations in technical units. Material science
analyzes the physical nature of fatigue and ways to determine the
number of cycles until a material fails under cyclic loading. In
this work, we will consider a compressor component of a gas
turbine subject to surface driven low-cycle fatigue (LCF). For
backgrounds on surface driven LCF failure mechanism with respect
to polycrystalline metal we refer to \cite{Harders_Roesler},
\cite{Vormwald}, \cite{Fedelich} and \cite{Sornette}.

\begin{figure}[t]
    \centering
\scalebox{0.28}{\input{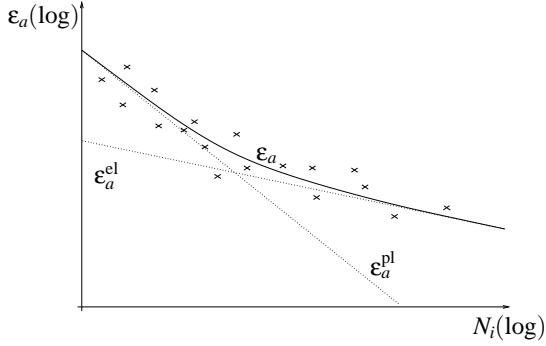}}
\caption{EN-DIAGRAM OF A STANDARDIZED SPECIMEN.}
    \label{Figure_woehler_diagram_strain}
\end{figure}

We now consider important methods of fatigue analysis and closely
follow Section 3 of \cite{Gottschalk_Schmitz} and
\cite{Harders_Roesler}. In fatigue specimen testing the number of
cycles until failure is determined. If the tests are strain
controlled so-called $E-N$ diagrams -- see Figure
\ref{Figure_woehler_diagram_strain} -- are created, where the
relationship between the strain amplitude $\varepsilon_a$ and
number $N_i$ of cycles until crack initiation is called W\"ohler
curve.
Usually, the range of cycles is subdivided into low-cycle fatigue
(LCF) and high-cycle fatigue (HCF). LCF loads are often strain
controlled, whereas HCF loads are mainly stress controlled so that
corresponding $S-N$ diagrams are analyzed.

For the purpose of analysis the strain amplitude $\varepsilon_a$
is subdivided into an elastic and plastic part where
$\varepsilon_a=\varepsilon_a^{el}+\varepsilon_a^{pl}$ holds. In
the LCF range the plastic part $\varepsilon_a^{pl}$ dominates,
whereas in the HCF range the elastic part
$\varepsilon_a^{\textrm{el}}$ plays a greater role. Introducing
the parameters fatigue strength $\sigma_f'$ and fatigue strength
exponent $b$ we present the so-called Basquin equation
$\varepsilon_a^{el}=\frac{\sigma_f'}{E}(2N_i)^b$ for the elastic
part, where $E$ is Young's modulus. Regarding the LCF range the
Coffin-Manson equation $\varepsilon_a^{pl}=\varepsilon_f'(2N_i)^c$
describes the denominating plastic part, where the parameters
$\varepsilon_i'$ and $c$ are called fatigue ductility and fatigue
ductility exponent, respectively. For a detailed discussion of the
physical origin of this equation we refer to \cite{Sornette}. The
combination of the previous equations leads to the
Coffin-Manson-Basquin (CMB) equation
\begin{equation}\label{materials_1.1}
\varepsilon_a=\varepsilon_a^{el}+\varepsilon_a^{pl}=\frac{\sigma_f'}{E}(2N_i)^b+\varepsilon_f'(2N_i)^c.
\end{equation}
The parameters can be calibrated according to the test data -- see
Figure \ref{Figure_woehler_diagram_strain} -- by means of maximum
likelihood methods, for example. Confer \cite{Georgii},
\cite{Escobar_Meeker}, \cite{Schott} and Subsection
\ref{subsection_FEA_results_calibration} below.

Structural design concepts with respect to LCF often consider the
component's surface position of highest stress and then analyze
the W\"ohler curve which corresponds to the conditions at that
surface position. Mostly safety factors are imposed to account for
the stochastic nature of fatigue, size effects\footnote{Note that
different geometries of test specimens lead to different W\"ohler
curves, confer \cite{Vormwald}.} and for uncertainties in the
stress and temperature fields. Note that sometimes several surface
positions of highest stress are considered which depends on the
component. This concept is called safe-life approach to fatigue
design and is often used in engineering as well as very similar
methods, confer \cite{Harders_Roesler}.

Note that several extensions exist to the CMB equation
(\ref{materials_1.1}). In particular, when approaching the HCF
region, mean stress effects are of increasing importance. The
modified Morrow equation is one approach that would consider such
effects. For further discussions confer \cite{Vormwald}.
%Having outlined the theoretical foundations for describing
%polycrystalline metal under cyclic loading we present the local
%and probabilistic model for LCF according to
%\cite{Gottschalk_Schmitz}.

\subsection{From Reliability Statistics to Probabilistic
LCF}\label{section_stochastics}

%As mentioned before the stochastic nature of fatigue requires
%several fatigue tests to obtain statistically significant
%statements for crack initiation times. An overview over fatigue
%tests and statistical methods can be found in \cite{Schott} and
%\cite{Vormwald}.
Now, we introduce the local and probabilistic model for LCF as
presented in \cite{Gottschalk_Schmitz} and
\cite{Siemens_Juelich_Wuppertal_Lugano} which can be derived from
reliability statistic, confer \cite{Escobar_Meeker}.

We model failure-time processes on continuous scale, although time
in our context is a number of load cycles and thereby an integer
number. Let $N$ denote a continuous random variable
%on some probability space with real non-negative values
which represents the time of crack initiation here identified with
failure of a system or component. If $P$ denotes the underlying
probability measure $F_N(n)=P(N\leq n)$ is the cumulative
distribution function and $f_N(n)=dF_N(n)/dn$ the density
function. The survival function is defined by
$S_N(n)=P(N>n)=1-F_N(n)$ and the hazard function by
\begin{equation*}
h(n)=\lim_{\Delta n\rightarrow0}\frac{P(n<N\leq n+\Delta
n|N>n)}{\Delta n}=\frac{f_N(n)}{1-F_N(n)},
\end{equation*}
confer \cite{Escobar_Meeker}. $h$ is also called hazard rate or
instantaneous failure rate function. For a small step $\Delta n$
the expression $h(n)\cdot\Delta n$ is an approximation for the
propensity of an object or system to fail in the next time step
$\Delta n$, given survival to time $n$. For a large number $m(n)$
of items in operation at time $n$ the product $m(n)\cdot h(n)$ is
approximately the number of failures per unit time. Defining the
cumulative hazard function $H(n)=\int_0^nh(t)dt$ one can show that
the survival function satisfies $S_N(n)=1-F_N(n)=\exp(-H(n))$.
This shows that a model ansatz for the hazard function leads to a
corresponding distribution function $F_N$.

%Let the failure-time be given by the random variable $N$.
Now, we introduce the crucial assumption for the local and
probabilistic model for LCF of\cite{Gottschalk_Schmitz} in the
case of polycrystalline metal. Consider LCF failure mechanism on
the component which is represented by the domain $\Omega$. We
assume that the surface zone that is affected from the crack
initiation process of a single LCF crack is small with respect to
the surface of the component. This surface zone corresponds to
faces of a few grains. As long range order phenomena are unusual
in polycrystalline metal, we pass to the following assumption:

%That is why we call the previous assumption the property of
%spatial additivity of hazard rates.

\noindent \textbf{Assumption (L)}\\
In any surface region $A\subseteq\partial\Omega$ the corresponding
hazard rate $h_A$ is a local functional of the elastic
displacement field $u$ in that particular region with
\begin{equation}
h_{A}(n)=\int_{A} \rho(n;\nabla \textbf{u},\nabla^2\textbf{u})\,
dA.
\end{equation}
Here, $\nabla \textbf{u}$ is the Jacobian matrix and
$\nabla^2\textbf{u}$ the Hessian of $\textbf{u}$. Note that the
loads that we consider in the given context (mainly elastic
plastic stresses and strains) can all be expressed as functions of
$\nabla \textbf{u}$. The integrand $\rho$ is called hazard density
function.

Assumption (L) is also called the property of spatial additivity
of hazard rates. A motivation of (L) and more detailed backgrounds
can be found in \cite{Gottschalk_Schmitz}. For mathematical
simplicity, we restrict ourselves to the case where only the
dependence on elastic strains (or equivalently stresses) is taken
into account.
%Modeling damage times in the spirit
%of elastic support factors depending on $\nabla^2u$ via
%$\chi^*=\frac{|\nabla \sigma_v|}{\sigma_v}$ might be an
%interesting option for the future.

In the case of inhomogeneous strain fields assumption (L) implies
$h(n)=\int_{\partial\Omega}\rho(n;\varepsilon^e)\,dA$
for some hazard density function $\rho$. %The specific form of
%$\rho$ depends on the failure mechanism and will be discussed in
%the following.
Because of
$F_N(n)=1-\exp(-H(n))=1-\exp\left(-\int_0^nh(t)\,dt\right)$ we
obtain for the probability of failure in $\partial\Omega$ until
cycle $n$:
\begin{equation}\label{materials_2.1}
F_N(n)=1-\exp\left(-\int_0^n\int_{\partial\Omega}\rho(t;\varepsilon^e)\,dAdt\right).
\end{equation}
The ansatz (\ref{materials_2.1}) can also be derived by the
Poisson point process with $\rho(n;\varepsilon^e(x))$ as the
intensity measure. For details on point processes we refer to
\cite{Baddeley} and \cite{Klenke}. The advantage of the point
process is that also the probability of a given number of cracks
initiations in $A \subseteq
\partial\Omega$ within $n$ load cycles can be computed via
\begin{equation}
P({\rm number~of~crack~ initiations~on~}A=q)=e^{-z}\frac{z^q}{q!}
\end{equation}
for $z=\int_0^n\int_{A}\rho(t;\varepsilon^e)\,dAdt$. But if cracks
have grown sufficiently large they will mutually influence their
local stress fields and thus the approach will break down.

We now establish a link to deterministic LCF analysis via CMB
equation (\ref{materials_1.1}) which leads to an appropriate
choice for the hazard density function $\rho$. We assume that the
number $N$ of cycles to crack initiation are Weibull distributed
which can be realized by the choice of a Weibull hazard ansatz
\begin{equation}\label{materials_3.1}
%\begin{split}
\rho(n;\mathbf{x})=\rho(n;\varepsilon^e(\mathbf{x}))=\frac{m}{N_{\textrm{det}}(\varepsilon^e(\mathbf{x}))}\left(\frac{n}{N_{\textrm{det}}(\varepsilon^e(\mathbf{x}))}\right)^{m-1}.
%\end{split}
\end{equation}
Here, $m$ is the Weibull shape and
$N_{\textrm{det}}(\varepsilon^e)$ the Weibull scale parameter
which is supposed to depend on the elastic strain tensor
$\varepsilon^e(\mathbf{x})$of the BVP
(\ref{Thermoelasticity.1.0}). Combining (\ref{materials_2.1}) and
(\ref{materials_3.1}) leads to the following model of
\cite{Gottschalk_Schmitz}:

\textbf{(Local and Probabilistic Model for LCF)}\\
Let the scale field
$N_{\textrm{det}}(\mathbf{x})=N_{\textrm{det}}(\varepsilon_a(\mathbf{x}))$,
$\mathbf{x}\in\partial\Omega$, be the solution of the CMB equation
(\ref{materials_1.1})
\begin{equation}\label{materials_5.1}
\varepsilon_a(\mathbf{x})=\frac{\sigma_f'}{E}(2N_{\textrm{det}}(\mathbf{x}))^{b}+\varepsilon_f'(2N_{\textrm{det}}(\mathbf{x}))^{c},
\end{equation}
where $\varepsilon_a(\mathbf{x})$ is computed from
$\varepsilon^e(\mathbf{x})$ via\footnote{Confer Sections
\ref{section_thermoelasticity_plasticity} and
\ref{section_fatigue}.} linear isotropic elasticity, from the von
Mises stress $\sigma_v(\mathbf{x})$, from
$\sigma_a(\mathbf{x})=SD^{-1}(\sigma_v(\mathbf{x})/2)$ according
to Neuber shakedown and from the Ramberg-Osgood equation with
$\varepsilon_a(\mathbf{x})=RO(\sigma_a(\mathbf{x}))$.
%calculated from $\varepsilon^e(x)$ via the stress-strain relation,
%calculation of the elastic von Misies stress $\sigma_v^e$, the
%shake down and the Ramberg-Osgood formula.
Then, the local and probabilistic model for LCF is given by the
cumulative distribution function
\begin{equation}\label{materials_4.1}
F_N(n)=1-\exp\left(-\int_0^n\int_{\partial\Omega}\frac{m}{N_{\textrm{det}}}\left(\frac{s}{N_{\textrm{det}}}\right)^{m-1}dAds\right).
\end{equation}
for $n\geq0$ and some $m\geq1$, which yields the probability for
LCF crack initiation in the interval $[0,n]$.

The shape parameter $m$ determines the scatter of the distribution
where small values for $m\geq1$ correspond\footnote{$0<m\leq1$ is
not realistic for fatigue.} to a large scatter and where the limit
$m\to\infty$ is the deterministic limit. Note that the Weibull
hazard function can be easily replaced by any other differentiable
hazard function with scale parameter $N_{\textrm{det}}$.

The CMB parameters of the model are not the same as obtained from
fitting standard specimen data. We calibrate them by means of
usual maximum likelihood methods, confer Section
\ref{section_application}, \cite{Georgii} and
\cite{Escobar_Meeker}. Furthermore, note that volume driven
fatigue could be considered as well by replacing the surface
integral in (\ref{materials_4.1}) with a volume integral whose
integrand only differs by different material parameters. For a
discussion of volume driven fatigue such as HCF confer
\cite{Harders_Roesler}.

With respect to materials engineering, the local and probabilistic
model for LCF has significant advantages compared to the safe-life
approach to fatigue design: The model bypasses the standard
specimen approach and takes size effects into account, i.e.
results from arbitrary geometries under LCF failure mechanism can
be employed to calibrate our model and every position of the
surface of an engineering part is considered by a surface integral
which does not need information on W\"ohler curves of a specific
specimen. Thereby inhomogeneous stress fields are taken into
account.

%%%%%%%%%%%%%%%%%%%%%%%%%%%%%%%%%%%%%%%%%%%%%%%%%%%%%%%%%%%%%%%%%%%%%%
\section{FINITE ELEMENT ANALYSIS AND POSTPROCESSING}\label{section_Postprocessor}

In this section we describe our finite-element postprocessor which
computes the distribution function with respect to fatigue life.
First, we briefly introduce into concepts of FEA. Then, we discuss
numerical integration and backgrounds on the postprocessor.

\subsection{Finite Element Analysis and Lagrange Elements}\label{subsection_FEA}

In order to numerically solve the BVP (\ref{Thermoelasticity.1.0})
of linear elasticity on a three-dimensional polyhedron $\Omega$ we
apply FEA, where a so-called weak formulation is considered on a
finite-dimensional space of functions, confer
\cite{Finite_Elements_Ern}, \cite{Braess} and
\cite{Finite_Elements_Ciarlet_1}. This function space and the
geometry $\Omega$ are described by a mesh of finite elements,
where each element consists of a three-dimensional compact and
connected set $T$ and of a finite set of basis functions
$\Pi=\{\psi_1,\dots,\psi_{n_{sh}}\}$ associated to a set of nodes
$\{\textbf{a}_1,\dots,\textbf{a}_{n_{sh}}\}$ in $T$. In the
following $\{T,\Pi\}$ denotes a finite element and the functions
of $\Pi$ are called shape functions. The FEA solution of the BVP
(\ref{Thermoelasticity.1.0}) restricted to $T$ is then a certain
linear combination of the shape functions.
%Under further assumptions on $\Pi$ and $\Sigma$ there is basis
%$\{\psi_1,\dots,\psi_{n_{sh}}\}$ in $\Pi$ where every $\psi_i$ can
%be uniquely assigned to one side condition of $\Sigma$. Each
%$\psi_i$ is called shape function.

Very popular finite elements are the Lagrange finite elements
which are contained in most FEA packages. In our example case of a
compressor disk we will use Abaqus 6.9-2 and Lagrange elements of
Serendipity class (C3D20) where the shape functions are of the
form
\begin{equation}\label{ICS_Serendipity}
\psi(\mathbf{x})=\sum_{\stackrel{0\leq i_1,i_2,i_3\leq
2}{i_1+i_2+i_3\leq 3}}\alpha_{i_1,i_2,i_3}\cdot x_1^{i_1}x_2^{i_2}
x_3^{i_3}%\right|\alpha_{i_1,i_2,i_3}\,\textrm{real}.
\end{equation}
for coefficients $\alpha_{i_1,i_2,i_3}$ determined by
$\psi_i(\mathbf{a}_j)=\delta_{ij}$.
%The side conditions $\Sigma$ are given by nodes
%$\textbf{a}_1,\dots,\textbf{a}_{n_{sh}}\in T$ such that
%$p(\mathbf{a}_1),\dots,p(\mathbf{a}_{n_{sh}})$ has prescribed
%values for all functions $p$ in $\Pi$.

Now, principles of mesh and finite element generation are briefly
introduced, where we closely follow \cite{Finite_Elements_Ern}.

\textbf{(Mesh)}\\
The mesh of a domain $\Omega$ is given by the union of compact and
connected sets $\{T_m\}_{1\leq m\leq N_{el}}$ with
$\overline{\Omega}=\bigcup_{m=1}^{N_{el}}\,T_m$ and all interiors
nonempty and pairwise disjoint. $h_T$ denoting the greatest
distance of two points in $T$ we define
$\mathcal{K}_h=\{T_m\}_{1\leq m\leq N_{el}}$ where
$h=\max\left\{h_T\,|\,T\in\{T_m\}_{1\leq m\leq N_{el}}\right\}$.

The starting point for the construction of meshes is a geometric
reference cell $\hat{T}$ from which the different sets are
generated.

\textbf{(Generated Mesh, Reference Cell)}\\
A mesh $\mathcal{K}_h$ of a domain $\Omega$ is called generated if
there exists a diffeomorphism\footnote{A diffeomorphisms is
continuously differentiable map whose inverse exists and is
continuously differentiable as well.} $\Upsilon_T$ for every
$T\in\mathcal{K}_h$ and a fixed set $\hat{T}$ such that
$\Upsilon_T(\hat{T})=T$. The set $\hat{T}$ is called reference
cell.

Now, we define the so-called geometric transformation which can
provide the previous diffeomorphism, confer
\cite{Finite_Elements_Ern}:

\textbf{(Geometric Transformation)}\\
Let $\{\hat{T},\hat{\Pi}\}$ be a finite element with shape
functions $\{\hat{\psi}_1,...,\hat{\psi}_{n_{sh}}\}$ and let
$\{\mathbf{a}_1,...,\mathbf{a}_{n_{sh}}\}$ be three-dimensional
nodes. Then the geometric transformation denotes the map
\begin{equation}\label{ICS_geometric_transformation}
%\Upsilon:\hat{T}\longrightarrow\mathbf{R}^d,\quad\quad
\hat{\mathbf{x}}\longmapsto\Upsilon(\hat{\mathbf{x}})=\sum_{i=1}^{n_{sh}}\mathbf{a}_i\,\hat{\psi}_i(\hat{\mathbf{x}}).
\end{equation}

In most cases a mesh generator produces a list of
$n_{\textrm{geo}}$ times $N_{el}$ nodes such as
%\begin{equation}
$\{\mathbf{a}_1^m,...,\mathbf{a}_{n_{\textrm{geo}}}^m\}_{1\leq
m\leq N_{el}}$.
%\end{equation}
Here, $N_{el}$ is the number of elements again. The vectors
$\{\mathbf{a}_1^m,...,\mathbf{a}_{n_{\textrm{geo}}}^m\}$ are also
defined as the geometric nodes of the $m$-th element. If
$\Upsilon_m$ is the geometric transformation corresponding to the
reference cell $\hat{T}$ and to
$\{\mathbf{a}_1^m,...,\mathbf{a}_{n_{\textrm{geo}}}^m\}$ with
$m=1,\dots,N_{el}$ we set $T_m=\Upsilon_m(\hat{T})$. So we obtain
$m$ elements and can write down the following definition:

\textbf{(Geometric Reference Finite Element)}\\
If a finite element $\{\hat{T},\hat{\Pi}_{\textrm{geo}}\}$ is used
for mesh generation by means of given nodes
\begin{equation}
\{\mathbf{a}_1^m,...,\mathbf{a}_{n_{\textrm{geo}}}^m\}_{1\leq
m\leq N_{el}},%\,\subset\,\mathbf{R}^d,
\end{equation}
its geometric transformations $\Upsilon_m$ and setting
$T_m=\Upsilon_m(\hat{T})$ for $m=1,\dots,N_{el},$ it is called the
geometric reference finite element.

\textbf{(Generation of Lagrange Elements)}\\
Finally, we consider the generation of Lagrange elements
$\{T,\Pi\}$ from a reference Lagrange element
$\{\hat{T},\hat{\Pi}\}$. Here, $\{\hat{\mathbf{a}}_i\}_{1\leq
i\leq n_{sh}}$ denotes the Lagrange nodes of
$\{\hat{T},\hat{\Pi}\}$ and
$\{\hat{\psi}_1,...,\hat{\psi}_{n_{sh}}\}$ the corresponding shape
functions which are assigned to the nodes via
$\hat{\psi}_i(\hat{\mathbf{a}}_j)=\delta_{ij}$. In this case
$\hat{\Pi}=\hat{\Pi}_{\textrm{geo}}$ and given nodes
$\{\mathbf{a}_1^m,...,\mathbf{a}_{n_{\textrm{geo}}}^m\}_{1\leq
m\leq N_{el}}$ provide a mesh $\mathcal{K}_h$ with reference cell
$\hat{T}$, geometric transformations $\Upsilon_T$ and sets
$T=\Upsilon_T(\hat{T})$. The nodes of $T$ satisfy
$\mathbf{a}_i=\Upsilon_T(\hat{\mathbf{a}}_i)$ for
$i=1,\dots,n_{sh}$. In order to construct $\Pi$ consider for
functions $v$ the map
\begin{equation}\label{ICS_example}
%\phi_T:V(T)\longrightarrow V(\hat{T}),\quad\quad
v\longmapsto\phi_T(v)=v\circ \Upsilon_T
\end{equation}
which is linear and invertible. Then, the shape functions are
given by $\Pi=\{\phi_T^{-1}(\hat{p})\,|\,\hat{p}\in\hat{\Pi}\}$.
%The side conditions $\Sigma$ are given by nodes
%$\textbf{a}_1,\dots,\textbf{a}_{n_{sh}}$ such that
%$v(\mathbf{a}_1),\dots,v(\mathbf{a}_{n_{sh}})$ has prescribed
%values for all functions $v$ in $\Pi$.

%The subsequently constructed vector spaces lead to approximation
%spaces such as $V_h$ of (\ref{FEA.2.0}).

%As we employ Lagrange elements of the form (\ref{ICS_Serendipity})
%in our example case we consider finite elements of degree $2$
%according to \cite{Finite_Elements_Ern}. Now, additionally assume
%that the boundary of zero displacements $\partial\Omega_D$ is the
%union of mesh faces. Then, Proposition 3.82 in
%\cite{Finite_Elements_Ern} implies that the approximation error
%$\textbf{u}-\textbf{u}_h$ converges to zero with respect to
%appropriate integral norms if the mesh size $h$ tends to zero.
%Here, $\textbf{u}$ is the solution of the weak formulation of BVP
%(\ref{Thermoelasticity.1.0}) and $\textbf{u}_h$ the numerical FEA
%solution. If $\textbf{u}$ is more regular\footnote{Here, more
%regularity means $u\in [H^{l}(\Omega)]^3\cap V$ for $l>1$ in terms
%of Sobolev spaces. For Sobolev spaces and corresponding integral
%norms confer Appendix B.3 in \cite{Finite_Elements_Ern}.} better
%convergence rates can be realized. A good approximation
%$\textbf{u}_h$ will be obtained if $\Omega$ can be subdivided into
%finite elements with sufficiently small mesh size $h$.

\subsection{Postprocessing}\label{subsection_postprocessor}

Recall that the local and probabilistic model for LCF is given by
the cumulative distribution function (\ref{materials_4.1}).
%\begin{equation}\label{ICS_recall}
%F_N(n)=1-\exp\left(-\int_0^n\int_{\partial\Omega}\frac{m}{N_{det}}\left(\frac{s}{N_{det}}\right)^{m-1}dAds\right).
%\end{equation}
%Here, $N_{\textrm{det}}(\varepsilon^e)$ is the scale parameter and
%is computed as described in the remarks after
%(\ref{materials_3.1}).
In the following we address numerical integration and explain how
the linearized strain field
$\varepsilon^e(\mathbf{x}),\mathbf{x}\in\Omega,$ is computed by
means of linear elastic results
from FEA with Lagrange elements. %In our example case of a
%compressor disk Lagrange elements of Serendipity class are
%employed.

The surface integral in (\ref{materials_4.1}) is numerically
computed by quadrature formulae. Note that we can also numerically
compute volume integrals, where in (\ref{materials_4.1})
$\partial\Omega$ is replaced by $\Omega$. This can be of higher
interest if we extent our model to volume driven failure mechanism
such as HCF. Because locations of stress concentrations will have
a major contribution to our probability functional we have to
consider higher nonlinearities in the integrand by our numerical
approach. Therefore, we will use quadrature formulae of higher
order.
%Furthermore, note that often the method of reduced
%integration is used in FEA packages where certain integration
%points are neglected to reduce needed computing performance.
%%This results in the need for local evaluations of the linearized strain
%%rate tensor $\varepsilon^e(x)$, where $x\in\Omega$ is arbitrary.

\begin{figure}[t]
    \centering
\scalebox{0.28}{\input{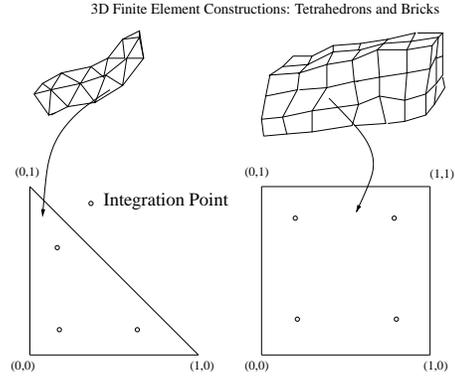}}
\caption{GEOMETRIC TRANSFORMATION AND INTEGRATION POINTS ON THE
UNIT TRIANGLE AND RECTANGLE.}
    \label{Figure_integration_points}
\end{figure}

Since $\Omega$ is a polyhedron which consists of finite elements
we can conduct the integration of (\ref{materials_4.1}) on every
face of a finite element which contributes to the boundary of
$\partial\Omega$ and then sum up each of that integral values.
Moreover, the geometric transformation
(\ref{ICS_geometric_transformation}) can be used as a chart to
transform each face integration on the corresponding map area
which is the unit triangle in case of tetrahedrons or the unit
rectangle in case of bricks. Also see Figure
\ref{Figure_integration_points}, where exemplary geometric
transformations and integration points are depicted. For more
details on surface integration consider the comments to
(\ref{ICS_surface_integration}) below and the term of integration
over submanifolds as described in \cite{Forster_III}, for example.
Substitution for multiple variables is employed in case of volume
integrals and thus the element integrals are transformed to the
unit tetrahedron or unit brick. %For substitution for multiple
%variables confer \cite{Forster_III}.

\textbf{(Quadrature of Order $k$)}\\
Let $K$ be an integration domain, let the $l_q$ numbers
$\omega_1,\dots,\omega_{l_q}$ be denoted as weights and the $l_q$
points $\mathbf{\xi}_1,\dots,\mathbf{\xi}_{l_q}\in K$ as
integration points (quadrature points). The weights and
integration points are called quadrature of order $k$, if $k$ is
the largest integer such that
\begin{align}
\begin{split}
&\int_K p(\mathbf{x})dx=\sum_{l=1}^{l_q}\omega_l
p(\mathbf{\xi}_l)\quad\textrm{for
all}\\
&p\in\mathbf{P}_k=\left\{\left.\sum_{\stackrel{0\leq
i_1,\dots,i_d\leq k}{i_1+\dots+i_d\leq
k}}\alpha_{i_1,\dots,i_d}x_1^{i_1}\dots
x_d^{i_d}\right|\alpha_{i_1,\dots,i_d}\,\textrm{real}\right\}.
\end{split}
\end{align}
With respect to intervals $[a,b]$ Table
\ref{table_quadrature_rectangle} shows corresponding weights and
integration points. Confer Table 8.1. and 8.2. in
\cite{Finite_Elements_Ern} as well, where Table 8.2 contains
quadrature formulae for triangles. Note that a quadrature on an
interval can be used to obtain a quadrature on a rectangle
$[a,b]\times[c,d]$ by subdividing the multidimensional integral
into one-dimensional integrals. This results into the four
integration points of the unit rectangle in Figure
\ref{Figure_integration_points}. %, where we used two integration
%points for the $x_1$- and $x_2$-direction.
Similarly we obtain quadrature formulae for bricks and thus can
numerically integrate over volumes of finite elements. For volume
quadratures on tetrahedrons and quadratures of higher orders we
again refer to Section 8 in \cite{Finite_Elements_Ern}.

\begin{table}[t]
\caption{QUADRATURE ON THE INTERVAL $K=[a,b]$ WITH
$\tilde{m}=(a+b)/2$ AND $\tilde{\delta}=b-a$.}
\begin{center}
\label{table_quadrature_rectangle}
\begin{tabular}{c l l}
& & \\ % put some space after the caption
\hline
$k_q\,\,\,\, l_q$ \vline & Int. Points $\xi_l$ & Weights $\omega_l$ \\
\hline
$1\quad1$  \vline & $\tilde{m}$ & $\tilde{\delta}$ \\
\hline
$3\quad2$  \vline & $\tilde{m}\pm\frac{\tilde{\delta}}{2\sqrt{3}}$ & $\frac{1}{2}\tilde{\delta}$ \\
\hline
$5\quad3$  \vline & $\tilde{m}\pm\frac{\tilde{\delta}}{2}\sqrt{\frac{3}{5}}$ & $\frac{5}{18}\tilde{\delta}$ \\
$\quad\quad$  \vline & $\tilde{m}$ & $\frac{8}{18}\tilde{\delta}$ \\
\hline
$7\quad4$  \vline & $\tilde{m}\pm\frac{\tilde{\delta}}{2}\left(\sqrt{(15+2\sqrt{30})/35}\right)$ & $\left(\frac{1}{4}-\frac{1}{12}\sqrt{\frac{5}{6}}\right)\tilde{\delta}$ \\
$\quad\quad$  \vline & $\tilde{m}\pm\frac{\tilde{\delta}}{2}\left(\sqrt{(15-2\sqrt{30})/35}\right)$ & $\left(\frac{1}{4}+\frac{1}{12}\sqrt{\frac{5}{6}}\right)\tilde{\delta}$ \\
\hline
\end{tabular}
\end{center}
\end{table}

%\begin{table}[t]
%\caption{QUADRATURE ON A TRIANGLE OF AREA $S$.}
%\begin{center}
%\label{table_quadrature_triangle}
%\begin{tabular}{c l l l}
%& & \\ % put some space after the caption
%\hline
%$k_q\,\,\,\, l_q$ \vline & Int. Points (Barycentric) & Multiplicity & Weights $\omega_l$ \\
%\hline
%$1\quad1$  \vline & $\left(\frac{1}{3},\frac{1}{3},\frac{1}{3}\right)$ & $1$ & $S$ \\
%\hline
%$2\quad3$  \vline & $\left(\frac{1}{2},\frac{1}{2},0\right)$ & $3$ & $\frac{1}{3}S$ \\
%\hline
%$3\quad4$  \vline & $\left(\frac{1}{3},\frac{1}{3},\frac{1}{3}\right)$ & $1$ & $-\frac{9}{16}S$ \\
%$\quad\quad$  \vline & $\left(\frac{1}{5},\frac{1}{5},\frac{3}{5}\right)$ & $3$ & $\frac{25}{48}S$ \\
%\hline

%$4\quad6$  \vline & $\left(a_i,a_i,1-2a_i\right),i=1,2$ & $3$ & $b_iS,i=1,2$\\
%$\quad\quad$  \vline & $a_1\approx0.446$ &  & $b_1\approx0.223$\\
%$\quad\quad$  \vline & $a_2\approx0.092$ &  & $b_2\approx0.110$\\
%\hline
%\end{tabular}
%\end{center}
%\end{table}

%As previously mentioned we employ quadratures of higher order than
%those in FEA packages like Abaqus and Ansys. Moreover, we want to
%have the flexible choice to extend our tool with new quadratures
%and new finite elements. Therefore, we cannot use the integration
%points and field values at these points of an FEA package and must
%evaluate the field at any possible point in every finite element.
%Thus, we have to employ the functional part $\Pi$ of each finite
%element.

Considering the finite element generation, we use the coordinates
of the nodes, the displacements of the nodes and the connectivity
(information on which nodes belong to a specific element).
Moreover, we have to employ the explicit form of the shape
functions $\{\hat{\psi}_0,\dots,\hat{\psi}_{n_{sh}}\}$ of the
reference Lagrange element $\{\hat{T},\hat{\Pi}\}$. Then, we
obtain the geometric transformation for every finite element
$\{T_l,\Pi_l\}$ according to (\ref{ICS_geometric_transformation})
%\begin{equation}\label{ICS_geometric_trafo}
%\Upsilon_{T_l}:\hat{T}\longrightarrow T_l,\quad\quad
%\hat{x}\longmapsto\Upsilon_{T_l}(\hat{\mathbf{x}})=\sum_{k=1}^{n_{sh}}\textbf{a}_k\,\hat{\psi}_k(\hat{\mathbf{x}}),
%\end{equation}
and the local displacements on every $T_l$ according to
(\ref{ICS_example}):
\begin{equation}\label{ICS_displacement}
\hat{\mathbf{x}}\longmapsto
(\textbf{u}_{T_l}\circ\Upsilon_{T_l})(\hat{\mathbf{x}})=\sum_{k=1}^{n_{sh}}\textbf{u}_k\,\hat{\psi}_k(\hat{\mathbf{x}}),
\end{equation}
with $\textbf{u}_k$ displacement vectors at each node of
$\{T_l,\Pi_l\}$.

Now, the computation of the linearized strain rate tensor
$\varepsilon^e$ is explained. We will express the fields in
dependence of $\hat{\mathbf{x}}\in\hat{T}$. Let $\mathbf{x}$ be
the coordinates that describe the locations of $T_l$. Because of
the geometric transformation we have
$\mathbf{x}=\Upsilon_{T_l}(\hat{\mathbf{x}})$ for some
$\hat{\mathbf{x}}\in\hat{T}$. Recall the form of the linearized
strain rate tensor:
\begin{equation*}
\varepsilon_{ij}(\mathbf{x})=\frac{1}{2}\left(\frac{\partial
u_i}{\partial x_j}(\mathbf{x})+\frac{\partial u_j}{\partial
x_i}(\mathbf{x})\right),\quad i,j\in\{1,2,3\},
\end{equation*}
with $u_i$ components of displacement vector \textbf{u}. In order
to express the strain completely in coordinates $\hat{\mathbf{x}}$
of the reference element $\{\hat{T},\hat{\Pi}\}$, one has to
consider the chain rule and to implement
\begin{equation}\label{ICS_displacement_reference_derivative}
\frac{\partial\left(\textbf{u}_{T_l}\circ\Upsilon_{T_l}\right)_i}{\partial
\hat{x}_j}(\hat{\mathbf{x}})
=\sum_{k=1}^{n_{sh}}\left(\textbf{u}_k\right)_i\,\frac{\partial\hat{\psi}_k}{\partial\hat{x}_j}(\hat{\mathbf{x}}),\quad
i\in\{1,2,3\},
\end{equation}
(which follows from (\ref{ICS_displacement})) and
\begin{equation}\label{ICS_geometric_transformation_derivative}
\frac{\partial\left(\Upsilon_{T_l}\right)_i}{\partial
\hat{x}_j}(\hat{\mathbf{x}})
=\sum_{k=1}^{n_{sh}}\left(\textbf{a}_k\right)_i\,\frac{\partial\hat{\psi}_k}{\partial\hat{x}_j}(\hat{\mathbf{x}}),\quad
i,j\in\{1,2,3\}.
\end{equation}
(\ref{ICS_geometric_transformation_derivative}) leads to the
Jacobian $\nabla\Upsilon_{T_l}(\hat{x})$ of the geometric
transformation $\Upsilon_{T_l}$ of
(\ref{ICS_geometric_transformation}).

%For $i\in\{1,2,3\}$ the chain rule leads to
%\begin{align}\label{ICS_displacement_derivative}
%\begin{split}
%&\left(\frac{\partial\left(\textbf{u}_{T_l}\right)_i}{\partial
%x_1}(\mathbf{x}),\frac{\partial\left(\textbf{u}_{T_l}\right)_i}{\partial
%x_2}(\mathbf{x}),
%\frac{\partial\left(\textbf{u}_{T_l}\right)_i}{\partial x_3}(\mathbf{x})\right)=\\
%&\left(\frac{\partial\left(\textbf{u}_{T_l}\circ\Upsilon_{T_l}\right)_i}{\partial
%\hat{x}_1}(\hat{\mathbf{x}}),
%\frac{\partial\left(\textbf{u}_{T_l}\circ\Upsilon_{T_l}\right)_i}{\partial
%\hat{x}_2}(\hat{\mathbf{x}}),
%\frac{\partial\left(\textbf{u}_{T_l}\circ\Upsilon_{T_l}\right)_i}{\partial
%\hat{x}_3}(\hat{\mathbf{x}})\right)\\
%&\cdot\left(D\Upsilon_{T_l}(\hat{\mathbf{x}})\right)^{-1}
%\end{split}
%\end{align}

Considering material parameters such as Young's modulus $E$ and
Poisson's ratio $\nu$, we can now compute with respect to the
coordinates $\hat{\mathbf{x}}\in\hat{T}$ the stress tensor, the
von Mises stress (\ref{Thermoelasticity.7.1}), the elastic-plastic
von Mises stress $\sigma_v$ according to Neuber shakedown
(\ref{Thermoelasticity.9.1}), the comparison strain
$\varepsilon_v$ according to the Ramberg-Osgood equation
(\ref{Thermoelasticity.8.1}) and finally $N_{\textrm{det}}$
according to the CMB approach (\ref{materials_5.1}).

Recalling (\ref{materials_4.1}), we explain the computation of
$F_N(n)$: Let $\mathcal{F}_{ij}$ be the $j$-th face of the $i$-th
element for $i=1,\dots,N_{el}$ and $j=1,\dots,N_F$, where $N_F$ is
the number of faces of the considered element
type. %\footnote{Bricks: $N_F=6$. Tetrahedrons: $N_F=4$.}
Let $\delta_{\mathcal{F}_{ij},\partial\Omega}$ be $1$ if
$\mathcal{F}_{ij}$ is a subset of the surface $\partial\Omega$
with area greater than $0$, otherwise let it be $0$. Then,
\begin{align}\label{ICS_integration}
\begin{split}
&\int_0^n\int_{\partial\Omega}\frac{m}{N_{\textrm{det}}}\left(\frac{s}{N_{\textrm{det}}}\right)^{m-1}dAds=n^m\int_{\partial\Omega}\left(\frac{1}{N_{\textrm{det}}}\right)^{m}dA\\
&=n^m\sum_{i=1}^{N_{el}}\sum_{j=1}^{N_F}\left[\delta_{\mathcal{F}_{ij},\partial\Omega}\cdot\int_{\mathcal{F}_i}\left(\frac{1}{N_{\textrm{det}}}\right)^{m}dA\right].
\end{split}
\end{align}
In case of bricks we integrate over the unit rectangle as map
area:
\begin{equation}\label{ICS_surface_integration}
\int_{\mathcal{F}_{ij}}\left(\frac{1}{N_{\textrm{det}}}\right)^{m}dA=\int_\Box\left(\frac{1}{N_{\textrm{det}}(\gamma_{ij}(\mathbf{s}))}\right)^{m}\sqrt{g^{\gamma_{ij}}(\mathbf{s})}ds
\end{equation}
with charts
$\gamma_{i1}(\mathbf{s})=\Upsilon_{T_i}(0,s_1,s_2),\gamma_{i2}(\mathbf{s})=\Upsilon_{T_i}(1,s_1,s_2),\gamma_{i3}(\mathbf{s})=\Upsilon_{T_i}(s_1,0,s_2),
\gamma_{i4}(\mathbf{s})=\Upsilon_{T_i}(s_1,1,s_2),\gamma_{i5}(\mathbf{s})=\Upsilon_{T_i}(s_1,s_2,0),\gamma_{i6}(\mathbf{s})=\Upsilon_{T_i}(s_1,s_2,1)$
for $\mathbf{s}=(s_1,s_2)\in[0,1]\times[0,1]$ and with the Gram
determinant
$g^{\gamma_{ij}}(\mathbf{s})=\det\left(\nabla\gamma_{ij}(\mathbf{s})^T\nabla\gamma_{ij}(\mathbf{s})\right)$,
where $\nabla\gamma_{ij}$ is the Jacobian of $\gamma_{ij}$. In
case of tetrahedrons we integrate over the unit triangle $\Delta$
as map area with charts
$\gamma_{i1}(\mathbf{s})=\Upsilon_{T_i}(0,s_1,s_2),\gamma_{i2}(\mathbf{s})=\Upsilon_{T_i}(s_1,0,s_2),\gamma_{i3}(\mathbf{s})=\Upsilon_{T_i}(s_1,s_2,0),
\gamma_{i4}(\mathbf{s})=\Upsilon_{T_i}(s_1,s_2,1-s_1-s_2)$. In
case of volume integral the Gram determinant
(\ref{ICS_integration}) has to be replace by
$|\det\left(D\Upsilon_{T_i}\right)|$ and the integration is
conducted over the whole reference brick and tetrahedron,
respectively.
%and bricks (\ref{ICS_integration}) becomes
%\begin{align}\label{ICS_integration_volume}
%\begin{split}
%&n^m\sum_{i=1}^{N_{el}}\left[\int_{T_i}\left(\frac{1}{N_{\textrm{det}}(\mathbf{x})}\right)^{m}dx\right]\\
%&=n^m\sum_{i=1}^{N_{el}}\left[\int_{\Box}\left(\frac{1}{N_{\textrm{det}}(\Upsilon_{T_i}(\hat{\mathbf{x}}))}\right)^{m}
%|\det\left(D\Upsilon_{T_i}(\hat{\mathbf{x}})\right)|d\hat{\mathbf{x}}\right].
%\end{split}
%\end{align}
%Here, $\int_{\Box}$ means the integration over the reference brick
%and will be replaced by $\int_{\Delta}$, the integration over the
%reference tetrahedron, if the elements are of tetrahedron type.

Now, we use quadrature formula $(\xi_l,\omega_l)_{l=1,\dots,l_q}$
to numerically compute the cumulative distribution function
(\ref{materials_4.1}). Considering (\ref{ICS_integration}) and
(\ref{ICS_surface_integration}) leads to:
\begin{align}\label{ICS_integration_approx}
\begin{split}
&\quad\,\, F_N(n)\\%=1-\exp\left(-\int_0^n\int_{\partial\Omega}\frac{m}{N_{\textrm{det}}}\left(\frac{s}{N_{\textrm{det}}}\right)^{m-1}dAds\right)\\
&\approx1-\exp\left(-n^m\sum_{i=1}^{N_{el}}\sum_{j=1}^{N_F}\delta_{\mathcal{F}_{ij},\partial\Omega}\sum_{l=1}^{l_q}
\frac{\sqrt{g^{\gamma_{ij}}(\xi_l)}}{\left(N_{\textrm{det}}(\gamma_{ij}(\xi_l))\right)^m}\cdot\omega_l\right)\\
&=1-\exp\left(-\left(\frac{n}{\eta}\right)^m\right),
\end{split}
\end{align}
where we defined the scale value
\begin{equation}\label{ICS_scale_variable}
    \eta=\left(\sum_{i=1}^{N_{el}}\sum_{j=1}^{N_F}\delta_{\mathcal{F}_{ij},\partial\Omega}\sum_{l=1}^{l_q}
\frac{\sqrt{g^{\gamma_{ij}}(\xi_l)}}{\left(N_{\textrm{det}}(\gamma_{ij}(\xi_l))\right)^m}\cdot\omega_l\right)^{-1/m}.
\end{equation}
In case of bricks use Table \ref{table_quadrature_rectangle} and
in case of tetrahedrons Table 8.2 in \cite{Finite_Elements_Ern}.
Equation (\ref{ICS_integration_approx}) shows that the Weibull
hazard ansatz in the local and probabilistic LCF model leads to a
Weibull distribution for the component $\Omega$ under surface
driven and strain-controlled LCF failure mechanism. Here, the
scale value $\eta$ is computed according to chosen finite element
types and quadrature formula. As mentioned in Subsection
\ref{section_stochastics} the parameters of the Ramberg-Osgood,
the CMB equation and the shape parameter $m$ have to be calibrated
via LCF-test results.

In this work the main input for FEA postprocessing are the
coordinates and displacements of the nodes and the connectivity.
Moreover, values for the parameters of our model must be given
which depend on the considered material. Note that we do not
employ additional information of FEA packages on which faces
$\mathcal{F}_{ij}$ contribute to the surface $\partial\Omega$ with
a non-vanishing
area, %On the one hand this is more convenient as our method is
%thereby more independent of FEA packages. On the other hand
%Identifying the surface requires a greater implementation effort
%and more computing capacity.
so that we identify the surface only via the main FEA output. This
identification of the surface can numerically fail if finite
elements are significantly distorted. But in most applications
there are at most a few of such distorted elements and their
contribution to the overall probability of crack initiation can
often be neglected.

\section{LCF-LIFE ESTIMATION FOR A COMPRESSOR DISK}\label{section_application}

\begin{figure}[t]
    \centering
\scalebox{0.2}{\includegraphics{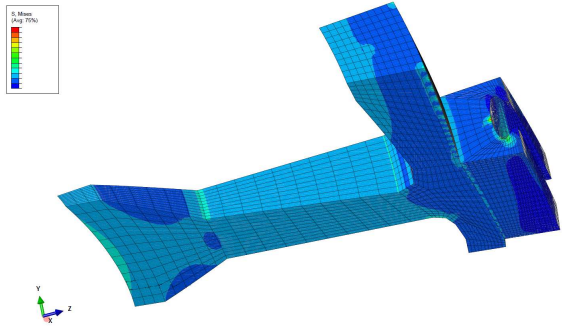}} \caption{FEA RESULTS
OF ABAQUS 6.9-2 FOR THE VON MISES STRESS FIELD OF THE COMPRESSOR
DISK.}
    \label{Figure_disk}
\end{figure}

In this section we consider a linear elastic Abaqus model of a
compressor disk, see Figure \ref{Figure_disk}, and estimate its
Weibull distribution with respect to LCF life by means of our
postprocessor.

\subsection{FEA Model and Parameter
Calibration}\label{subsection_FEA_results_calibration}

In the following we consider two states of the compressor disk:
The shutdown state and the operating state. Before the gas turbine
is activated the compressor disk is subjected to a homogeneous
temperature field but to no stresses. This constitutes the
shutdown state. In the operating state the disk is subjected to an
inhomogeneous stress field. The temperature field is still
homogeneous, but its value has increased. It is assumed that the
fields are stationary which is only an approximation regarding the
real operating state of the disk. As a conservative estimation the
temperature field of the shutdown state is set equal to that one
of the final state, otherwise thermal mechanical fatigue (TMF) has
to be considered. The transition from the shutdown state to the
operating state and then back to the shutdown state is considered
as one load cycle. It is assumed that the shutdown and operating
state stay the same during the cycles which is again only an
approximation.

An Abaqus FEA model has been created to predict the displacement
and von Mises stress field of the compressor disk in the operating
state. The model consists of approximately 10.000 Lagrange
elements of the Serendipity class (C3D20) which is given in
(\ref{ICS_Serendipity}). The compressor disk is subjected to a
homogeneous temperature field in the operating state which is
considered in the FEA model by corresponding values for Young's
modulus and Poisson's ratio. The shutdown state is already at our
disposal as the displacement and stress fields are zero and the
temperature field is set equal to the homogeneous one of the
operating state. Thus, the main input for our postprocessor are
the coordinates of the nodes, the displacements of the nodes and
the connectivity of the FEA model prediction which determine the
approximation of the displacement field in the operating state and
the geometry of the disk.

Finally, values for the material parameters of the LCF model must
be transferred. For this purpose we first calibrate the LCF model
with respect to the disk material according to standardized LCF
tests. %, where the cumulative distribution function $F_N(n)$ of
%(\ref{ICS_integration_approx}) can be effectively simplified as
%the surface is subjected to homogeneous stress and temperature
%fields. Thus, the surface integral is only a multiplication with
%the surface area between the gauge length.
Considering the cumulative distribution function $F_N(n)$ of
(\ref{ICS_integration_approx}) and defining
$\eta=\left(\int_{\partial\Omega}N_{\textrm{det}}^{-m}dA\right)^{-1/m}$,
we obtain for the corresponding density function $f_N(n)$ the
expression
\begin{equation}\label{results_1.0}
f_N(n)=\frac{d}{dn}F_N(n)=\frac{m}{\eta}\left(\frac{n}{\eta}\right)^{m-1}\exp\left[-\left(\frac{n}{\eta}\right)^m\right].
\end{equation}
We subsume the parameters of the model in a vector $\theta$, which
includes the parameters of Ramberg-Osgood and of CMB and the
Weibull shape parameter $m$. The experimental data set for $q$
strain-controlled LCF tests is given by
$\left\{n_i,\varepsilon_i,\partial\Omega_i\right\}_{i=1,\dots,q}$.
Here, $n_i$ is the number of cycles until crack initiation and
$\varepsilon_i$ is the strain on the gauge surface
$\partial\Omega_i$. We estimate (calibrate) $\theta$ by means of
maximum likelihood. The log-likelihood function is defined as
\begin{align}
\begin{split}
&\log\left(\mathcal{L}\left(\{(n_i,\varepsilon_i,\partial\Omega_i)\}_{i\in\{1,\dots,q\}}\right)[\theta]\right)\\
&=\sum_{i=1}^q\log(f_N(\partial\Omega_i,\varepsilon_i)(n_i)[\theta]).
\end{split}
\end{align}
Let $\hat{\theta}$ denote the likelihood estimator, then
$\hat{\theta}$ is given by
\begin{align}\label{results_2.0}
\begin{split}
&\log\left(\mathcal{L}\left(\{(n_i,\varepsilon_i,\partial\Omega_i)\}_{i\in\{1,\dots,q\}}\right)[\hat{\theta}]\right)\\
&=\max_{\theta}\left\{\log\left(\mathcal{L}\left(\{(n_i,\varepsilon_i,\partial\Omega_i)\}_{i\in\{1,\dots,q\}}\right)[\theta]\right)\right\}.
\end{split}
\end{align}
Having optimized (\ref{results_2.0}) and found the estimator
$\hat{\theta}$ we can now apply our postprocessor to the FEA model
of the compressor disk.

%In this example case we use the optimization software package
%DAKOTA 5.0 to find the estimator $\hat{\theta}$. Now, we can apply
%our postprocessor to the FEA model of the compressor disk.

\subsection{Results of the Probabilistic Approach to LCF}

\begin{figure}[t]
    \centering
\scalebox{0.55}{\includegraphics{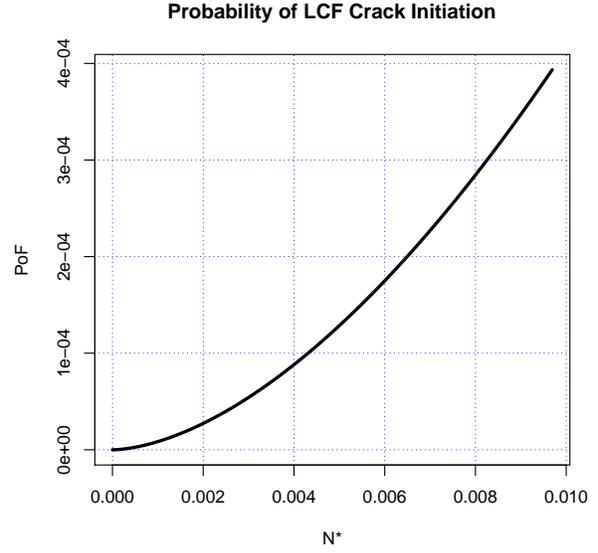}}
\caption{WEIBULL DISTRIBUTION FOR LCF CRACK INITIATION ON THE
COMPRESSOR DISK.}
    \label{Figure_PoF}
\end{figure}

The most important result of the numerical integration in our
model is the scale parameter $\eta$ of the Weibull distribution
(\ref{ICS_integration_approx}) which yields the probability for
LCF crack initiation until cycle $N$. The corresponding Weibull
shape parameter $m$ is already estimated by the calibration of the
previous section. The Weibull distribution is shown in Figure
\ref{Figure_PoF} for low numbers of cycles compared to $\eta$.
Here, failure is defined by the initiation of the first LCF crack
on the component and PoF denotes probability of failure and $N^*$
multiples of $\eta$.

\begin{figure}[t]
    \centering
\scalebox{0.14}{\includegraphics{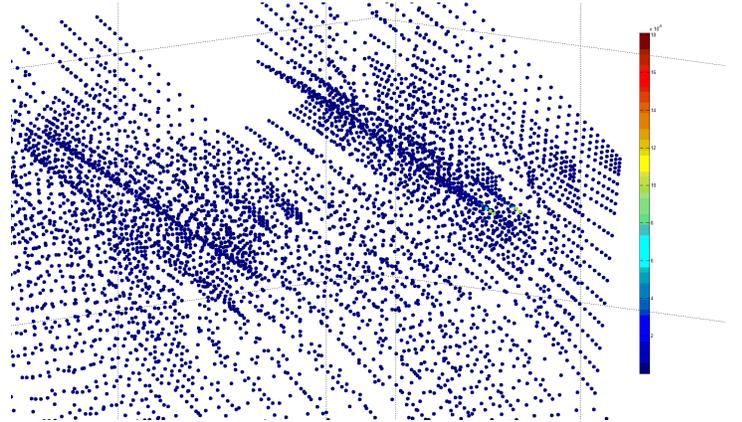}}
\caption{CRACK INITIATION DENSITY.}
    \label{Figure_intensity}
\end{figure}

%With respect to more design relevant cycles Figure
%\ref{Figure_PoF_design} shows the PoF in the range $N^*=0$ to
%$N^*=6.000$.
For $N^*=3.231\cdot10^{-3}$ the PoF is $6.142\cdot10^{-3}\%$, for % $N^*=0.003231$
example. If we consider that the compressor disk consists of 44
disk segments, where one is shown in Figure \ref{Figure_disk}, the
probability for the initiation of the first LCF crack on the
compressors disk is $0.270\%$. From a design perspective one
decides which PoF is acceptable and then chooses the corresponding
number of allowable shutdown and service cycles. Figure
\ref{Figure_intensity} shows the crack initiation density at $N^*$
on the disk's surface. This density corresponds to the local
expectation value for the number of crack initiations. One can see
that the crack initiation density is very much localized at the
bearing flank. As we first compute the Weibull scale for every
boundary face we can quantify this localization. In case of the
previous value for $N^*$, the $21$ faces with the greatest crack
initiation density have a combined PoF of already
$5.531\cdot10^{-3}\%$ which is more than $90\%$ of the PoF for all
boundary faces.

%\begin{figure}[t]
%    \centering
%\scalebox{0.55}{\includegraphics{PoF_design_2.eps}}
%\caption{WEIBULL DISTRIBUTION FOR LCF CRACK INITIATION -- DESIGN
%RANGE.}
%    \label{Figure_PoF_design}
%\end{figure}

With respect to field data of already operating compressor disks
no failure has occurred so far. Our postprocessor predicts for
that number of disks and service cycles a low PoF. Note that the
assumptions for the shutdown and final state as well as the FEA
model only approximate the real operating conditions of these
compressor disks. In particular, consider that the field data
consists of information from different gas turbines whose service
conditions can be very different. Moreover, uncertainties in the
model parameters $\theta$ -- recall the previous section and
confer the end of this section -- influence the real PoF of the
compressor disk. Nevertheless the model is able to predict a low
PoF.

\begin{figure}[t]
    \centering
\scalebox{0.55}{\includegraphics{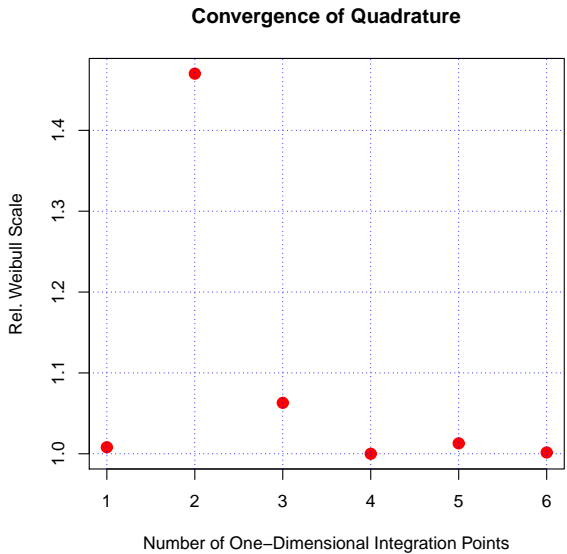}}
\caption{CONVERGENCE OF WEIBULL SCALE REGARDING NUMBER OF
INTEGRATION POINTS.}
    \label{Figure_quadrature_convergence}
\end{figure}

As a computational validation item of our tool we investigated
whether the numerical integration converges regarding the order of
the chosen quadrature. Since the FEA model consists of nonlinear
Lagrange elements of Serendipity class the postprocessor employs
quadratures of Table \ref{table_quadrature_rectangle}, confer
Subsection \ref{subsection_postprocessor}. We computed the Weibull
scale $\eta$ for every quadrature till the order of $k=11$, i.e.
till $l_q=6$ integration points in every dimension and thus $36$
integration points on the unit rectangle. Figure
\ref{Figure_quadrature_convergence} shows the results for each
Weibull scale divided by the estimated limit and depending on the
number $i=1,\dots,6$ of integration points in single dimension.
Using $16$ integration points on the unit rectangle, approximately
results in a converged value for $\eta$ which indicates a highly
nonlinear behavior of the integrand $N_{\textrm{det}}^{-m}$ of the
probability function (\ref{ICS_integration_approx}). This also
justifies the choice of quadratures of higher order.

The cumulative distribution function $F_N(n)$ in
(\ref{ICS_integration_approx}) depends on the parameters $\theta$
which are calibrated as described in Subsection
\ref{subsection_FEA_results_calibration}. Because this calibration
is a statistical estimation for $\theta$ depending on LCF test
data, there are uncertainties for the values of $\theta$. This
effects the real PoF but is not considered by our method so far.
Taking this additional uncertainty mathematically into account can
be realized by Bootstrap methods or Maximum Likelihood asymptotic
theory, for example. Confer \cite{Escobar_Meeker} for more details
on these methods.

Additional to the consideration of uncertainties in the model
parameters $\theta$, the extension of our model to inhomogeneous
temperature fields will be important for risk estimation for LCF
crack initiation on components such as turbine blades. For this
purpose a reliable temperature model for LCF is needed. Also note
that the model could be extended to consider HCF, TMF and
non-stationary FEA. From a design perspective there is also the
interesting possibility to optimize the PoF
(\ref{ICS_integration_approx}) with respect to the shape $\Omega$,
i.e. find a design $\Omega$ under certain constraints such that
the surface integral and PoF (\ref{ICS_integration_approx}) is
minimized. This is also called optimal reliability, confer
\cite{Gottschalk_Schmitz}. Because the integrand is sufficiently
regular under additional smoothness assumptions there is even a
link to gradient-based shape optimization, confer
\cite{Shape_Optimization_Haslinger} and \cite{Sokolowski_Zolesio}.
This could accelerate computational optimization efforts
significantly.

% Here's where you specify the bibliography style file.
% The full file name for the bibliography style file
% used for an ASME paper is asmems4.bst.
\bibliographystyle{asmems4}

%%%%%%%%%%%%%%%%%%%%%%%%%%%%%%%%%%%%%%%%%%%%%%%%%%%%%%%%%%%%%%%%%%%%%%
\begin{acknowledgment}
This work has been supported by the German federal ministry of
economic affairs BMWi via an AG Turbo grant. We wish to thank the
gasturbine technology department of the Siemens AG for stimulating
discussions and many helpful suggestions.
\end{acknowledgment}

%%%%%%%%%%%%%%%%%%%%%%%%%%%%%%%%%%%%%%%%%%%%%%%%%%%%%%%%%%%%%%%%%%%%%%
% The bibliography is stored in an external database file
% in the BibTeX format (file_name.bib).  The bibliography is
% created by the following command and it will appear in this
% position in the document. You may, of course, create your
% own bibliography by using thebibliography environment as in
%

% Here's where you specify the bibliography database file.
% The full file name of the bibliography database for this
% article is asme2e.bib. The name for your database is up
% to you.
%\bibliography{asme2e}

%%%%%%%%%%%%%%%%%%%%%%%%%%%%%%%%%%%%%%%%%%%%%%%%%%%%%%%%%%%%%%%%%%%%%%

\end{document}